# On the length of Golomb Ruler: A function construction approach based on difference triangle[*]

Yanqing Wang, Xiaoming Li

School of Management, Harbin Institute of Technology
Harbin 150001, China
yanqing@hit.edu.cn

**Abstract** Since the significance of Golomb Ruler Problem in some context, we proposes a function construction approach based on difference triangle to generate near-optimal Golomb rulers. Let $x_1, x_2, \ldots x_n$ be an increasing sequence of integers, where $x_1 = 0$, which satisfies the following conditions: if $|x_i - x_j| = |x_p - x_q|$ then $\{i, j\} = \{p, q\}$. Our objective is to find the order of minimum $x_n$ for any given $n$. In this paper, the two results in a paper are both improved. In addition, it will be shown that the length of Golomb Ruler have been shortened to a half, and that the satisfying sequence can not be generated by such a quadratic formula as $x_i = ai^2+bni+ci+dn^2+en+f$ for any rational $a, b, c, d, e$ and $f$.

## 1. Introduction

The concept of Golomb Rulers is a variation of graceful graphs. A graph of $e$ edges is said to be graceful if its nodes can be numbered by the integers from $\{0, 1, \ldots e\}$ such that the induced numbers on edges, i.e., the positive differences of numbers assigned to their end nodes, are all different. Lots of attempt has been made on determining of whether a graph is graceful or not. It has been proved that every complete graph of $n > 4$ nodes is not graceful[2]. But unknown are more than known it is not even known if all tree graphs are graceful[3].

If we do not put any restriction on the ruler's length and only require that all distances between pairs of marks be different, the shortest ruler, for a given $n$, is called Golomb Ruler.

The idea of Golomb Ruler has a great variety of applications[4][5]. For example, a Golomb Ruler may correspond to optimal antennas geometry[6]. Unfortunately, we haven't found an efficient algorithm to generate Golomb Rulers. We do not even know how long it should be in terms of order. In the following, we will shorten the length of the Golomb Ruler in [7], in other words, we have found a shorten Golomb Ruler.

---

[*] This paper is an English version of the reference [1] by and large. Some minor mistakes have been corrected here. In addition, the improvements to this paper are being undertaken recently.

For our purpose, let's rephrase the problem of finding a Golomb Ruler, i.e. determining $n$ numbers $0 = x_1 < x_2 < \ldots < x_n$ such that:

(1) If $|x_i - x_j| = |x_p - x_q|$ then $\{i, j\} = \{p, q\}$

(2) $x_n$ is minimized.

For convenience, we call sequence $x_1, x_2, \ldots x_n$ with $0 = x_1 < x_2 < \ldots < x_n$ a graceful sequence if it satisfies the condition (1). It is clear that any graceful sequence must have $x_n \geqslant C_n^2$ since there are $C_n^2$ induced differences [7] (in fact, $x_n > C_n^2$ for $n > 4$ [2]). And it isn't difficult to verify that $x_i = 2^{i-1} - 1$ generates graceful sequences for any $n$ [7], but unfortunately $x_n$ might be too big. In what follows, we'll present a formula, which generates a graceful sequence for any $n$ and $x_n$ is of order $n^3$. Compare the length with that of [7], and you should find that length has been shortened to a half. Also we shall give a proof that such a formula as $x_i = ai^2 + bni + ci + dn^2 + en + f$ cannot serve our purpose and we conjecture that $O(x_n) > n^2$.

## 2. Definitions and a lemma

### 2.1 Definitions

Definition 1. Integer sequence $0 = x_1 < x_2 < \ldots < x_n$ is called a graceful sequence if the differences of any pair of elements are unique. i.e. if $|x_i - x_j| = |x_p - x_q|$ then $\{i, j\} = \{p, q\}$.

Definition 2. Following [7], to any graceful sequence of integers, the differences between elements of $\{x_i\}$ can be arranged in a so-called Difference Triangle (DT), where an element $t_{i,j}$ of DT is defined by $t_{i,j} = x_{i+1} - x_{i+1-j}$, and it is in the $i^{th}$ row, the $j^{th}$ column, here $1 \leqslant j \leqslant i \leqslant n - 1$. Thus, the graceful sequence consists of $x_1 = 0$ and all the elements in the diagonal of corresponding DT.

### 2.2 Lemma

**Lemma 1**. Let $a = pN + r$, $b = qN + s$, for any positive number $p, q, r, s$ and N, if $N > r$, $N > s$ and $r \neq s$, then when N is given, $a \neq b$ must be true.

Proof: As $a = pN + r$, $b = qN + s$ thus $a - b = (p - q) N + (r - s)$. Because $r < N$, $s < N$, it is clear

that $|r - s| < N$.

Case 1, $p = q$. Here $a - b = r - s$. $r \neq s$, so $a \neq b$ must be true.

Case 2, $p \neq q$. Supposed $a = b$, then $|p - q| N = |r - s|$, the left side, $p \neq q$ so $|p - q| \geq 1$, $|p - q| N \geq N$; while in the right side, $r < N$, $s < N$, so $|r - s| < N$ must be true, a contradiction.

<div align="right">Q. E. D.</div>

## 3 Theorems

**Theorem 1**. The length of Golomb Rulers of $n$ marks is bounded by $\dfrac{(n-1)((n-1)^2 + 1)}{2}$.

Proof: Set $x_i = \dfrac{(i-1)(i-2)n}{2} + i - 1$, ($i = 1, 2, \ldots n$). Our fist step is to construct a DT. It's not hard to understand if any two elements in DT are different then our construction is successful. First, look at the following interesting Difference Triangle $DT_1$ together:

```
0n+1
0n+1   0n+2
2n+1   3n+2    3n+3
3n+1   5n+2    6n+3    6n+4
4n+1   7n+2    9n+3    10n+4   10n+5
5n+1   9n+2    12n+3   14n+4   15n+5   15n+6
6n+1   11n+2   15n+3   18n+4   20n+5   21n+6   21n+7
 …      …       …       …       …       …       …       …
```

It's easy to find that the graceful sequence, which generates the above $DT_1$, consists of 0 and all the elements in the diagonal of $DT_1$.

The elements in every column are increasing from top to end, so it's clear that there can't be two same elements in one column. There are $n - 1$ columns totally, and any element in the $i^{th}$ column mod $n$ equals to $i$ ($i = 1, 2, 3, \ldots n - 1$). So if we prove that any pair of elements in two different columns can't be equal, we can say that $DT_1$ satisfied our condition.

Now select any pair of elements from two different columns as you wish. Suppose one of them is $a = pn + j1$, the other is $b = qn + j2$, here $j1 \neq j2$. Because there are $n - 1$ rows and $n - 1$ columns totally, $j1 < n$ and $j2 < n$. Based on Lemma 1, $a \neq b$ is always correct. It proves that the above $DT_1$ is a graceful one. The upper boundary of above Difference Triangle $DT_1$ is $x_n = \dfrac{(n-1)(n-2)n}{2} + n - 1 =$

$$\frac{(n-1)((n-1)^2+1)}{2}.$$

In above $DT_1$, if we replace $n$ with $n-2$, the conclusion is still true. The result is the same as the Theorem 1 in [7], it's $x_n = (n-1)(1 + \frac{(n-2)^2}{2})$.

Q. E. D.

**Theorem 2**. The length of Golomb Ruler with $n$ marks can be shortened to about a half. When $n$ is an odd integer, $x_n \leqslant (n-1) + \frac{(n-1)^2(n-2)}{4}$, else $x_n \leqslant (n-1) + \frac{n(n-1)(n-2)}{4}$.

Proof: In above $DT_1$, if we replace $n$ with N, our objective now is to decrease N as much as possible. Thus we set $x_i = C_{i-1}^2 N + i - 1$ ($i = 2, 3 \ldots n$). If $n$ is an odd number then set $N = \frac{n-1}{2}$, else if $n$ is an even one then $N = \frac{n}{2}$. So the updated Difference Triangle $DT_2$ is below:

```
1
N+1      N+2
2N+1     3N+2    3N+3
3N+1     5N+2    6N+3   6N+4
4N+1     7N+2    9N+3   10N+4   10N+5
5N+1     9N+2    12N+3  14N+4   15N+5   15N+6
 …        …       …      …       …       …       …      …
```

Now that N is $\frac{n}{2}$ or $\frac{n-1}{2}$, from the 1st column to the $(N-1)^{th}$ one, the elements in every column mod N are from 1 to N–1 correspondingly, and those in the $N^{th}$ column mod N are 0. From the $(N+1)^{th}$ column to $(n-1)^{th}$ column, when the elements mod N their remainders are from 1 to N–1 (or to zero when $n$ is odd). Therefore, the $DT_2$ is split into two big blocks (see the dotted line in $DT_2$). The first N columns form the first block; the other block is from the $(N+1)^{th}$ column to the $(n-1)^{th}$ one.

Based on Theorem 1, it's easy to find that, inside any block, no two elements are duplicated. According to Lemma 1, if two columns mod N have different remainders, all elements of these two columns are unique. So our next target is to prove that any two elements in the different blocks, which have the same remainders mod by N, can't be the same.

In above DT$_2$, when $j$ is given, $t_{i,j} = x_{i+1} - x_{i+1-j}$ will be increasing with $i$. Thus the largest element in the $j^{th}$ column is:

$$t_{n-1,j} = x_n - x_{n-j}$$
$$= [\frac{(n-1)(n-2)}{2} - \frac{(n-j-1)(n-j-2)}{2}] \cdot N + (n-1) - (n-j-1)$$
$$= [j(n-2) - \frac{j(j-1)}{2}] \cdot N + j;$$

While the smallest element of $(N+j)^{th}$ is:

$$t_{N+j,N+j} = C^2_{N+j} N + N + j = [\frac{(N+j)(N+j-1)}{2}] \cdot N + j;$$

If we can prove $t_{n-1,j}$ is always less than $t_{N+j,N+j}$ to any given $j$, our goal is achieved. That means that to any $j$ ($j = 1, 2, \ldots N$) the following formula must be true.

$$\frac{(N+j)(N+j-1)}{2} > j(n-2) - \frac{j(j-1)}{2} \qquad (*)$$

**Case 1**, $n$ is an odd. $N = \frac{n-1}{2}$, we want to prove

$$[\frac{n-1}{2} + j] \cdot [\frac{n-1}{2} + j - 1]/2 + 1 > j(n-2) - \frac{j(j-1)}{2}$$

$$2j^2 - (n-1)j + \frac{(n-1)(n-3)}{4} + 2 > 0$$

The parabola $y = 2j^2 - (n-1)j + \frac{(n-1)(n-3)}{4} + 2$ opens upwards, and $(n-1)^2 - [2(n-1)(n-3) + 16] = -(n-3)^2 - 12 < 0$, so $y$ is always more than zero. Our goal is accomplished.

**Case 2**, $n$ is an even. $N = \frac{n}{2} > \frac{n-1}{2}$. If we replace N with $\frac{n}{2}$ in formula (*), it is obviously that the formula is correct.

Now that DT$_2$ is constructed successfully, the upper boundary of Golomb Ruler is:

$$x_n \leqslant (n-1) + \frac{(n-1)^2(n-2)}{4}, \text{ when } n \text{ is an odd};$$

$$x_n \leqslant (n-1) + \frac{n(n-1)(n-2)}{4}, \text{ when } n \text{ is an even}.$$

Q. E. D.

**Theorem 3**. To any rational number $a, b, c, d, e$ and $f$, graceful sequence cannot be generated by such a quadratic polynomial $x_i = ai^2 + bni + ci + dn^2 + en + f$.

Proof: We are discussing Difference Triangle, in which $t_{i,j} = x_{i+1} - x_{i+1-j}$, $1 \leqslant j \leqslant i \leqslant n-1$, so

if we set $x_i = a(i-1)^2+bn(i-1)+c(i-1)$, ($i =1, 2, 3, \ldots, n$), generality will not be lost. Since that $x_i = a(i-1)^2+bn(i-1)+c(i-1)$, the corresponding $DT_3$ is shown below:

| | | | | |
|---|---|---|---|---|
| $1a+(bn+c)$ | | | | |
| $3a+(bn+c)$ | $4a+2(bn+c)$ | | | |
| $5a+(bn+c)$ | $8a+2(bn+c)$ | $9a+3(bn+c)$ | | |
| $7a+(bn+c)$ | $12a+2(bn+c)$ | $15a+3(bn+c)$ | $16a+4(bn+c)$ | |
| $9a+(bn+c)$ | $16a+2(bn+c)$ | $21a+3(bn+c)$ | $24a+4(bn+c)$ | $25a+5(bn+c)$ |
| … | … | … | … | … |

If we can prove that there must be two equal elements in $DT_3$ when $n$ becomes big enough, our conclusion is made.

First of all, we should exclude some conditions unnecessary to be considered.

If $a = 0$. All elements in a single column are all the same, so $a \neq 0$;

If $b = 0$. Based on the Theorem 2 in [7], it is obviously that sequence $\{x_i\}$ is not graceful.

Every element in DT must be more than zero. Let's look at $t_{1,1} = a + bn + c$. If $b < 0$, no matter how big $a$ and $c$ are, when $n$ becomes big enough, $t_{1,1}$ might be less than zero. That is unreasonable, so $b < 0$ can not be true.

In addition, from $t_{1,1} = a + bn + c > 0$, we can get $c > -a - 2b$ (when $n = 2$); from $t_{n-1,1} = (2n - 3)a + bn + c > 0$, we can get $(2a+b)n > 3a - c$. No matter how small $3a - c$ is, $2a + b > 0$ must be true.

On all account, we just need to prove that to any integer $a, b, c, a \neq 0, b > 0, c > -a - 2b, 2a + b > 0$, the sequence $\{x_i\}$ corresponding with $DT_3$ is NOT graceful.

Let's pick up two special elements $t_{i1,j1}$ and $t_{i2,j2}$ in $DT_3$, in which $i_1 = n - 1, j_1 = b + 1, i_2 = j_2 = 2a + b + 1$. If $a > 0$, $j_1$ is smaller column; otherwise $j_1$ is the bigger one.

$\because 2a^2+b^2+2ab+2b+(-a-2b) = (a+b)^2+ a^2-a \geqslant 0$

$\therefore (-a-2b) \geqslant -(2a^2+b^2+2ab+2b)$

$\therefore c > -a-2b \geqslant -(2a^2+b^2+2ab+2b)$

$\therefore (2a^2+b^2+2ab+2b) + c > 0$

So, we can find that

$2a^2+b^2+2ab+2b+(2a+b+1)+1+c-1 > 2a+b+1$. If $n= 2a^2+b^2+2ab+2a+3b+2+c$, then $n-1 > 2a+b+1 = j_2$;

Also, $c > -a-2b$, $n-1=c+2a^2+b^2+2ab+2a+3b+2-1 > -a-2b+2a^2+b^2+2ab+2a+3b+1=(a+b)^2+a^2+a+b+1 \geqslant b+1 = j_1$.

We can conclude that: If $a > 0$, then $n - 1 > j_2 > j_1 > 1$; If $a \leq 0$, then $n - 1 > j_1 \geq j_2 > 1$. So we can see that $j_2, j_1$ all exist.

Now let's see what will happen when $n = 2a^2+b^2+2ab+2a+3b+2+c$.

$$
\begin{aligned}
t_{i1,j1} - t_{i2,j2} &= (x_{i1+1} - x_{i1+1-j1}) - (x_{i2+1} - x_{i2+1-j2}) \\
&= ai_1^2 + (bn+c)i_1 - [ai_2^2 + (bn+c)i_2] \\
&\quad - [a(i_1-j_1)^2 + (bn+c)(i_1-j_1)] + a(i_2-j_2)^2 + (bn+c)(i_2-j_2) \\
&= a[i_1^2 - i_2^2 - (i_1-j_1)^2] + (j_1-j_2)(bn+c) + x_1 \\
&= a[i_1^2 - i_2^2 - (i_1-j_1)^2] - 2a(bn+c) \\
&= a[(n-1)^2 - (2a+b+1)^2 - (n-b-2)^2] - 2a(bn+c) \\
&= a[n^2 - 2n + 1 - 4a^2 - 4ab - 4a - b^2 - 2b - 1 - n^2 + 2bn - b^2 + 4n - 4b - 4] - 2a(bn+c) \\
&= -2a[2a^2 + b^2 + 2ab + 2a + 3b + 2 + c - n] = 0
\end{aligned}
$$

So when $n$ becomes big enough, $t_{i1,j1}$ and $t_{i2,j2}$ must be equal.

Q. E. D.

## 4. Conclusion

In this paper, some progresses are made on finding near-optimal Golomb rulers. A shorten ruler of order $O(x_n) \leq \dfrac{n^3}{4}$ is found. Moreover, a theorem prevents people from attempting in vain to find such a quadratic polynomial as $x_i = ai^2 + bni + ci + dn^2 + en + f$ to mark Golomb Ruler. However, this problem has not been solved completely, quite a lot of work should be done. For example, machine proving might be an assistant way to help us construct DT efficiently.

## References

[1] Yanqing Wang, Xiaoming Li. Some New Development on the Length of Golomb's Rulers. Journal of Harbin Institute of Technology, 1993, vol.25, no.4, pp: 30-34 (in Chinese).

[2] Golomb S W. How to Number a Graph. Graph Theory and Computing. Academic Press, New York, 1972:23-37.

[3] Harary F. Lectures in Graph Theory. Stevens Inst of Tech, Spring, 1985.

[4] Ayari N, Luong V, Jemai A. A hybrid genetic algorithm for Golomb ruler problem. The 2010 IEEE/ACS International Conference on Computer Systems and Applications (AICCSA), IEEE, 2010:

1-4.

[5] Bloom G S, Golomb S W. Numbered Complete Graphs, Unusual Rulers, and Assorted Applications. Proceedings of the International Conference on the Theory and Applications of Graphs, Spring – Verlag, 1976.

[6] Unnikrishna Pillai S, Yeheskel Bar-Ness Fred Haber. A New Approach to Aray Geometry for Improved Spatial Spectrum Estimation. Proceedings of the IEEE, 1985, 73(10).

[7] Xiaoming Li. On the Length of Optimal Non-redundancy Linear Arrays. Proc. of ICCS'90, Nov., 1990, Singapore.